\documentclass[12pt]{amsart}
\usepackage{a4wide}
\usepackage{amsmath,amsfonts,amssymb}
\usepackage{amscd}
\usepackage{latexsym}

\newcommand{\LL}{\ell\ell}
\newcommand{\N}{\mathbb{N}}

\newtheorem{theorem}{Theorem}[section]

\begin{document}
\pagestyle{myheadings}
\title[Coincidences between character values for symmetric groups]{
Coincidences between characters to hooks and 2-part partitions
on families arising from 2-regular classes
}
\author{Christine Bessenrodt}
\address{
Leibniz Universit\"at Hannover,
Institute for Algebra, Number Theory and Discrete Mathematics,
Welfengarten 1,
30167 Hannover, Germany
}
\email{bessen@math.uni-hannover.de}
\date{October 23, 2016.\\
2010 {\em Mathematics Subject Classification.}
Primary 20C30, Secondary 05E10.\\
{\em Key words and phrases.}
Symmetric groups, characters.}

\maketitle

\begin{abstract}
Strongly refining results by Regev, Regev and Zeilberger,
we prove surprising coincidences between
characters to 2-part partitions of size $n$
and characters to hooks of size $n+2$
on two related families obtained by extending
2-regular conjugacy classes.
\end{abstract}

\section{Introduction}

In a recent paper,
Alon Regev, Amitai Regev and Doron Zeilberger
\cite{RRZ}
studied sums of squares of character values for the symmetric groups.
They compared on the one hand such a sum to characters where the labelling partitions
had a restricted number of rows, and on the other hand such a sum to
characters labelled by hooks; their values were considered on related
conjugacy classes with ``many'' fixed points, and
surprising coincidences between the sums of their squares were found.
In a subsequent paper \cite{RZ}, these results were greatly generalized to
further families of conjugacy classes.

We have to introduce some notation before stating the
details of these results and then present our strong refinements.
A partition of a number $n\in \N_0$
is an (unordered) multiset of nonnegative integers
summing to $n$,
often written as a weakly decreasing sequence and omitting trailing zeros;
the nonzero integers in the partition are its parts.
For $n=0$, the only partition is the empty partition, denoted $(0)$.
We also use exponential notation and write $i^m$ for $m$ parts~$i$
in a partition.
For a partition $\alpha$ and some positive integers $a,b,c,\ldots$,
we write $(\alpha,a,b,c,\ldots)$ for the partition having the parts of
$\alpha$ together with the parts $a,b,c,\ldots$.
It is wellknown that the irreducible complex characters
of the symmetric group $S_n$ are canonically labelled by the partitions of $n$;
we refer to \cite{JK} for background on this.
We denote the irreducible character to the partition $\lambda$ by $\chi^\lambda$.
The value of the character $\chi^\lambda$
on the conjugacy class of elements of cycle type $\mu$  is denoted by $\chi^\lambda(\mu)$. For $n=0$, we set $\chi^{(0)}((0))=1$.

Using the notation of \cite{RRZ}, we set
$$
\psi^{(2)}_n(\mu)= \sum_{j=0}^{\lfloor n/2 \rfloor} \chi^{(n-j,j)}(\mu)^2
$$
and
$$
\phi^{(2)}_n(\mu)= \sum_{j=1}^{n} \chi^{(j,1^{n-j})}(\mu)^2 \:.
$$
In \cite{RRZ}, they observed (and proved) the {\em remarkable identity}
$$
\psi^{(2)}_n(31^{n-3})= \frac 12 \phi^{(2)}_{n+2}(321^{n-3})
$$
and stated: {\em It may be interesting to find a `natural'
reason for this `coincidence'.}

In a subsequent paper \cite{RZ}, Amitai Regev and Doron Zeilberger
greatly generalized the identity above.
We quote their result which gives a wealth of surprising relations between characters to hooks and 2-part partitions.
\begin{theorem}\label{thm:RZ}
Let $\mu_0$ be a partition with odd parts $a_1\ge a_2 \ge \cdots \ge a_s\ge 3$
and the parts $2,2^2,\ldots,2^{t}$ for some $t\ge 0$.
Let $\mu_0'$ be the partition with the same
odd parts $a_1\ge a_2 \ge \cdots \ge a_s\ge 3$
and the part $2^{t+1}$.
Then for every $n\ge m=|\mu_0|$ we have
$$
\psi^{(2)}_n(\mu_0, 1^{n-m})= \frac 12 \phi^{(2)}_{n+2}(\mu_0', 1^{n-m})\:.
$$
\end{theorem}

Here, we will explain these results by showing that a strong refinement
of this relation is true: we have indeed coincidences between the character
values themselves not just the sum of squares of the values!
More precisely, we will show:
\begin{theorem}\label{thm:main}
Let $\alpha$ be a partition with odd parts only.
Let $\mu=(\alpha,2,2^2,\ldots,2^{t})$ for some $t\ge 0$, $n=|\mu|$.
Let $\mu'=(\alpha,2^{t+1})$.
Then for any $j \in \{0,1,\ldots, n-j\}$ we have
$$
\chi^{(n-j,j)}(\mu)= \chi^{(n+2-j,1^j)}(\mu')\:.
$$
\end{theorem}

\bigskip

Theorem~\ref{thm:RZ} is easily implied by this result.
We explain this using the notation in the theorem above.
As $\mu'$ is the cycle type for a class of odd permutations, we have
$$\chi^{(n+2-j,1^j)}(\mu')= - \chi^{(j+1,1^{n+1-j})}(\mu')\:.$$
Note that for odd $n=2k+1$ we have
$\chi^{(n+2-(k+1),1^{k+1})}(\mu') = \chi^{(k+2,1^{k+1})}(\mu')=0$.

Hence we obtain in any case, using Theorem~\ref{thm:main},
$$
\sum_{j=0}^{\lfloor n/2 \rfloor} \chi^{(n-j,j)}(\mu)^2
=
\sum_{j=0}^{\lfloor n/2 \rfloor} \chi^{(n+2-j,1^j)}(\mu')^2
=
\sum_{j=\lfloor n/2 \rfloor +1}^{n+1} \chi^{(n+2-j,1^j)}(\mu')^2
$$
and hence Theorem~\ref{thm:RZ} follows.

\section{Preliminaries}

We want to use an observation from the 2-modular representation theory of
the symmetric groups; we refer to the books by James and Kerber \cite{J,JK} and Olsson's Lecture Notes \cite{JO-LN} for more on the background.
Here, this will not imply that we work at characteristic~2,
but that we consider the restriction of the characters to 2-regular elements, i.e., those whose cycle type is a partition with odd parts only.

\smallskip

For two characters $\chi,\psi$ of $S_n$,
we write $\chi \equiv \psi$ if the
characters $\chi,\psi$ coincide on all 2-regular elements.

\smallskip

A crucial tool is a relation between characters to hooks and 2-part partitions on 2-regular elements.
We recall this observation here (see \cite[p.93]{J}).

Take $x\ge y\ge 0$. Then by the Littlewood-Richardson rule we have
$$
(\chi^{(x)} \times \chi^{(y)})\uparrow^{S_{x+y}}
=\chi^{(x+y)}+\chi^{(x+y-1,1)}+\ldots+\chi^{(x,y)}
$$
and
$$
(\chi^{(x)} \times \chi^{(1^y)})\uparrow^{S_{x+y}}
=\chi^{(x+1,1^{y-1})}+\chi^{(x,1^y)} \:.
$$
Note that for $y=0$ this is also correct by interpreting a character
to a non-partition like $(x+2,1^{-1})$ as zero.

Clearly, the characters $\chi^{(y)}$ and $\chi^{(1^y)}$
coincide on all 2-regular elements as these are even permutations.
Hence we deduce
$$
\chi^{(x+1,1^{y-1})}+\chi^{(x,1^y)}
\equiv
\chi^{(x+y)}+\chi^{(x+y-1,1)}+\ldots+\chi^{(x,y)}\:.
$$
Still using our convention from above and
$$
\chi^{(x+2,1^{y-2})}+\chi^{(x+1,1^{y-1})}
\equiv
\chi^{(x+y)}+\chi^{(x+y-1,1)}+\ldots+\chi^{(x+1,y-1)}
$$
we obtain the crucial relation
\begin{eqnarray}\label{eq:hooks2part}
\chi^{(x,y)}
\equiv
\chi^{(x,1^y)}-\chi^{(x+2,1^{y-2})}
\:.
\end{eqnarray}

\medskip

We need some further preliminaries on a variation of $\beta$-numbers;
for some of the following we refer to \cite{JO-LN} for more on $\beta$-numbers.

In general, any finite set $X\subset \N_0$ is a $\beta$-set.
For a partition $\lambda$ of length $\ell$, the first column hook lengths
$h_{i1}=\lambda_i+\ell-i$, $i=1,\ldots,\ell$,
form a $\beta$-set $X_\lambda$ for $\lambda$.
Conversely, to $X=\{x_1>x_2>\ldots > x_m\}$ we associate the partition
$$\pi(X)= (x_1-(m-1),x_2-(m-2),\ldots,x_{m-1}-1,x_m)$$
where trailing zeros are ignored.

For $r\in \N$, the $r$-shift of a $\beta$-set $X$ is  $X^{+r}=\{x+r\mid x\in X\}\cup \{r-1,\ldots,0\}$.
Note that all shifts of $X$ have the same
associated partition; indeed, any partition corresponds to
a shift class of $\beta$-sets.

Considering $\beta$-sets instead of the partitions themselves
is useful for the process of removing hooks:
the removal of a hook of length $h$ from $\lambda$ corresponds
to subtracting $h$ from one of the numbers $x\ge h$ in its $\beta$-set
$X_\lambda$ such that $x-h\not\in X_\lambda$.

As we want to use $\beta$-numbers in the context of the Murnaghan-Nakayama formula,
we also want to keep track of the sign of the leg lengths of the hooks to be removed.
Instead of writing this out at each step, we will use a slight
variation of $\beta$-sets.
We will only use this here in the context of two-row partitions.

An ordered $\beta$-set $Z=((z_1,\ldots,z_m))$ is a sequence such that
the corresponding set is a $\beta$-set $X$ of size $m$;
the associated partition is $\pi(Z)=\pi(X)$.
The associated sign of the ordered $\beta$-set $Z$ is defined to be the sign of the permutation sorting the entries into decreasing order.
Then, when we start from the ordered version of $X_\lambda=\{x_1>x_2>\ldots > x_m\}$,
this being denoted by $Z_\lambda=((x_1,\ldots,x_m))$,
the subtraction of $h$ from $x_i$ (say) is recorded at position~$i$,
without reordering, and we obtain for the removal of the corresponding
hook $H$ from $\lambda$:
$Z=((x_1,\ldots,x_i-h,\ldots,x_m))$, with $\pi(Z)=\lambda/H$.
Denoting the leg length  of the hook $H$ by $\LL(H)$,
note that $(-1)^{\LL(H)}$ is the sign of $Z$.

We extend the notation for ordered $\beta$-sets a little further
when the sequences are of length~2,
and we define associated virtual characters.
We allow also formal expressions $((x,y))$ with $x$ or $y$ negative,
or with $x=y$; in these cases, we don't have associated partitions
and we define the associated virtual characters to be zero.
When the entries are different and nonnegative, we call the
ordered $\beta$-set proper.
In this case,
recalling how $\beta$-sets are associated with partitions, and keeping in mind that the ordered $\beta$-set also records a sign,
we are led to associate
virtual characters to ordered $\beta$-sets as follows:
$$
\text{the virtual character  to $((x,y))$ is } \quad
\chi^{((x,y))}= \left\{
\begin{array}{rl}
\chi^{(x-1,y)} & \text{if } x>y\ge 0\\
-\chi^{(y-1,x)} & \text{if } y>x\ge 0\\
\end{array}
\right. \:.
$$
By the relation observed in~(\ref{eq:hooks2part}),
we deduce for the restriction on 2-regular elements:
$$\chi^{((x,y))}\equiv
\left\{
\begin{array}{rl}
\chi^{(x-1,1^y)}-\chi^{(x+1,1^{y-2})}, & \text{if } x>y\ge 0\\
\chi^{(y+1,1^{x-2})}-\chi^{(y-1,1^x)}, & \text{if } y>x\ge 0\\
\end{array}
\right.
\:.
$$
Conjugating a partition labelling a character
gives a character with the same values on 2-regular elements,
thus in the second case
$$
\chi^{(y+1,1^{x-2})}-\chi^{(y-1,1^x)} \equiv \chi^{(x-1,1^{y})}-\chi^{(x+1,1^{y-2})} \:,
$$
and hence for any nonnegative integers $x\neq y$ we obtain
\begin{eqnarray}\label{eq:hooks2partgen}
\chi^{((x,y))}\equiv
\chi^{(x-1,1^y)}-\chi^{(x+1,1^{y-2})}
\:.
\end{eqnarray}
In fact, also for $x=y$ both sides are zero;
if $x$ or $y$ are negative, then also by our conventions both
sides are zero.
Hence relation (\ref{eq:hooks2partgen}) holds for all $x,y$.

\section{Proof of Theorem~\ref{thm:main}}

We recall the set-up of the theorem.
Let $\alpha$ be a partition with odd parts only,
let $\mu=(\alpha,2,2^2,\ldots,2^{t})$ for some $t\ge 0$,
and let $\mu'=(\alpha,2^{t+1})$.
Set $n=|\mu|$. Then for any $j \in \{0,1,\ldots, n-j\}$ we want
to show that
$$
\chi^{(n-j,j)}(\mu)= \chi^{(n+2-j,1^j)}(\mu')\:.
$$

Applying the Murnaghan-Nakayama formula $t$ times we have 
$$
\chi^{(n-j,j)}(\mu)=
\sum_{H_1} \sum_{H_2} \cdots \sum_{H_t}
(-1)^{\sum_{j=1}^t \LL(H_j)}\chi^{(n-j,j)\setminus H_t\cdots \setminus H_2 \setminus H_1} (\alpha)
$$
where the $j$-th sum 
runs over the hooks $H_j$ of length $2^j$ of the partition  \hbox{$(n-j,j)\setminus H_t\cdots \setminus H_{j+1}$}, for
$j\in \{1,\ldots,t\}$,
computing this from the innermost to the outermost sum, starting with hooks $H_t$  of length $2^t$ of $(n-j,j)$.

Transferring this into the language of ordered $\beta$-sets and keeping our conventions on the associated virtual characters in mind, this becomes
\begin{eqnarray}\label{eq:value}
\chi^{(n-j,j)}(\mu)=
\sum_{I\subseteq \{1,\ldots,t\}}
\chi^{((n+1-j-\sum_{i\in I} 2^i,j-\sum_{i\in I^c} 2^i))} (\alpha)
\end{eqnarray}
where $I^c=\{1,\ldots,t\}\setminus I$.
We note here that clearly the sign of every proper ordered $\beta$-set
$$((n+1-j-\sum_{i\in I} 2^i,j-\sum_{i\in I^c} 2^i))$$
is equal to the sign $(-1)^{\sum_{j=1}^t \LL(H_j)}$ coming from the removal
of the corresponding $t$ hooks.
Furthermore, we note that if $2^t,\ldots, 2^k$ have been subtracted from
$((n+1-j,j))$ and an ordered $\beta$-set $((x,x))$ is obtained
(which gives a zero character by definition), then continuing with further
subtractions will always produce
pairs of ordered $\beta$-sets $((x-a,x-b)), ((x-b,x-a))$ with cancelling pairs of
characters, so the total
contribution at the final level $t$ from such an intermediate term $((x,x))$  will
still be zero.

The $2^t$ ordered $\beta$-set labels appearing in the sum above are easily
determined, as any even number $m$ with $0 \le m \le 2^{t+1}-2$ can uniquely be
written as
$m=\sum_{i\in I} 2^i$ for a suitable subset $I\subseteq \{1,\ldots,t\}$.
So with $((u,v))=((n+1-j,j))$, in the sum above all labels
of the following form appear:
$$
((u-2k,v-(2^{t+1}-2-2k)))=((n+1-j-\sum_{i\in I} 2^i,j-\sum_{i\in I^c}2^i)),
0 \le k \le 2^t-1\:.
$$
Note that, of course, some of these labels will give a zero character
contribution.

Recalling that we want to evaluate the virtual characters appearing in the sum
in equation~(\ref{eq:value}) at the 2-regular element $\alpha$, we can then
apply the relations deduced in Section~2.
We recall relation~(\ref{eq:hooks2partgen}):
$$\chi^{((x,y))}\equiv
\chi^{(x-1,1^y)}-\chi^{(x+1,1^{y-2})}
\:.
$$
Note here, that when $x\in \{0,1\}$, then the positive term disappears,
and when $y\in \{0,1\}$, then the negative term disappears, i.e.,
$$
\chi^{((x,y))}\equiv
\left\{
\begin{array}{rl}
-\chi^{(x+1,1^{y-2})} & \text{if } x\in \{0,1\}\\
\chi^{(x-1,1^y)} & \text{if } y\in \{0,1\}
\end{array}
\right.
\:.
$$

We discuss first the generic situation where we assume that $u=n+1-j$ and $v=j$ are both of size at least $2^{t+1}-2$.
In this case we obtain for the virtual character in equation~(\ref{eq:value})
$$
\begin{array}{rcl}
\displaystyle
\sum_{I\subseteq \{1,\ldots,t\}}
\chi^{((n-j+1-\sum_{i\in I} 2^i,j-\sum_{i\in I^c} 2^i))}
&=&
\displaystyle
\sum_{k=0}^{2^t-1} \chi^{((u-(2^{t+1}-2-2k),v-2k))} \\[7pt]
&\equiv &
\displaystyle
\sum_{k=0}^{2^t-1} (\chi^{(u-(2^{t+1}-2-2k)-1,1^{v-2k})}-\chi^{(u-(2^{t+1}-2-2k)+1,1^{v-2k-2})})
\\[7pt]
&=&
\displaystyle
\chi^{(u-(2^{t+1}-2)-1,1^{v})} - \chi^{(u+1,1^{v-(2^{t+1}-2)-2})}
\\[7pt]
&=&
\displaystyle
\chi^{(n+2-j-2^{t+1},1^{j})} - \chi^{(n+2-j,1^{j-(2^{t+1}-2)-2})}
\end{array}
$$
Remembering that $\chi^\rho=0$ if $\rho$ is not a partition,
we thus arrive at
$$
\begin{array}{rcl}
\chi^{(n-j,j)}(\mu)
&=&
\chi^{(n+2-j-2^{t+1},1^{j})}(\alpha)
- \chi^{(n+2-j,1^{j-2^{t+1}})}(\alpha)
\\[7pt]
&=&
\chi^{(n+2-j,1^j)}(\alpha,2^{t+1}) =\chi^{(n+2-j,1^j)}(\mu')
\:.
\end{array}
$$
Thus the claim is proved in this case.

Now assume that $n+1-j=u< 2^{t+1}-2$, say $u=\epsilon_u+2k_u$,
for some $0\le k_u<2^t-1$ and $\epsilon_u \in \{0,1\}$;
note that $u+v=n+1 \ge 2^{t+1}-2$, so then $v> 2^{t+1}-2-2k_u$.
Then the first possibly nonzero contribution in the sum comes
from the label $((\epsilon_u,v-(2^{t+1}-2-2k_u)))$,
for $k=2^t-1-k_u$. As we have seen above, for this label we only get the negative
hook contribution which cancels with the next term in the sum as before,
as long as such a term exists;
the critical case with no cancellation occurs only when
$v-(2^{t+1}-2-2k_u)=\epsilon_v\in \{0,1\}$.
Note that $v=j \leq n-j =u-1< 2^{t+1}-2$, so an analogous symmetric
argument for the last possibly nonzero contribution in the sum
tells us that in the noncritical cases or in the critical case
with $\epsilon_u=\epsilon_v$,
we have
$$
\begin{array}{rcl}
\chi^{(n-j,j)}(\mu)
=0=
\chi^{(n+2-j,1^j)}(\alpha,2^{t+1})
\end{array}
\:.
$$
Now we consider the critical case
$v-(2^{t+1}-2-2k_u)=\epsilon_v\neq \epsilon_u$.
Here, we have only one nonzero contribution at level $t$,
from the ordered $\beta$-set $((\epsilon_u,\epsilon_v))$,
namely, $\chi^{(0)}$ when the $\beta$-set is $((1,0))$, i.e.,
$u=n+1-j$ is odd, and $-\chi^{(0)}$ when it is $((0,1))$, i.e.,
$u=n+1-j$ is even.
Furthermore, in the critical case we have
$2^{t+1}-2 =u+v-1=n$;
thus $\alpha=\emptyset$ and we obtain
$$
\begin{array}{rcl}
\chi^{(n-j,j)}(\mu)=\chi^{(n-j,j)}(2,2^2,\ldots,2^t)
=(-1)^{n-j}=(-1)^j=
\chi^{(n+2-j,1^j)}(2^{t+1})
=\chi^{(n+2-j,1^j)}(\mu')\:.
\end{array}
$$

If $u \ge 2^{t+1}-2$, but $j=v< 2^{t+1}-2$, then the first term coming from the sum is still $\chi^{(n+2-j-2^{t+1},1^j)}$
as in the generic case, but arguing as above
there will not be a negative contribution at the end, so
$$
\begin{array}{rcl}
\chi^{(n-j,j)}(\mu)
=\chi^{(n-j-(2^{t+1}-2),1^{j})}(\alpha) =
\chi^{(n+2-j,1^j)}(\alpha,2^{t+1})
=\chi^{(n+2-j,1^j)}(\mu')\:.
\end{array}
$$
Hence we are now done in all cases.  $\Box$

\end{document}